\documentclass{amsart}
\usepackage{hyperref}
\usepackage{amsmath, amsthm, amsfonts, mathtools, tikz-cd, amssymb}
\usepackage{enumitem}

\makeatletter
\def\paragraph{\@startsection{paragraph}{4}%
  \z@\z@{-\fontdimen2\font}%
  {\normalfont\bfseries}}
\makeatother
\makeatletter
\@namedef{subjclassname@2020}{%
  \textup{2020} Mathematics Subject Classification}
\makeatother

\title[Bounded cohomology of classifying spaces for families]{Bounded cohomology of classifying spaces for families of subgroups}
\author{Kevin Li}
\address{School of Mathematical Sciences, University of Southampton, Southampton SO17 1BJ, United Kingdom}
\email{kevin.li@soton.ac.uk}
\date{May 19, 2021}

\theoremstyle{definition}
\newtheorem{defn}{Definition}[section]
\newtheorem{ex}[defn]{Example}
\theoremstyle{plain}
\newtheorem{thm}[defn]{Theorem}
\newtheorem{lem}[defn]{Lemma}
\newtheorem{prop}[defn]{Proposition}
\newtheorem{cor}[defn]{Corollary}
\theoremstyle{remark}
\newtheorem{rem}[defn]{Remark}

\newcommand{\IZ}{\ensuremath\mathbb{Z}}
\newcommand{\IQ}{\ensuremath\mathbb{Q}}
\newcommand{\IR}{\ensuremath\mathbb{R}}
\newcommand{\F}{\ensuremath\mathcal{F}}
\newcommand{\G}{\ensuremath\mathcal{G}}
\newcommand{\TR}{\ensuremath\mathcal{TR}}
\newcommand{\FIN}{\ensuremath\mathcal{FIN}}
\newcommand{\VCY}{\ensuremath\mathcal{VCY}}
\newcommand{\AME}{\ensuremath\mathcal{AME}}
\newcommand{\ALL}{\ensuremath\mathcal{ALL}}
\newcommand{\calH}{\ensuremath\mathcal{H}}
\newcommand{\OFG}{{\ensuremath\mathcal{O}_\mathcal{F}G}}
\newcommand{\OTRG}{\ensuremath\mathcal{O}_\mathcal{TR} G}
\newcommand{\EFG}{\ensuremath E_\mathcal{F}G}
\newcommand{\EGG}{\ensuremath E_\mathcal{G}G}
\newcommand{\Mod}{\ensuremath\text{-}\mathbf{Mod}}
\newcommand{\enum}{\rm{(\roman*)}}
\newcommand{\spann}[1]{{\ensuremath \langle{#1}\rangle}}
\newcommand{\cell}{\textup{cell}}
\newcommand{\linf}{\ensuremath \ell^\infty}
\DeclareMathOperator{\Hom}{Hom}
\DeclareMathOperator{\bHom}{bHom}
\DeclareMathOperator{\bMap}{bMap}
\DeclareMathOperator{\Ext}{Ext}
\DeclareMathOperator{\res}{res}
\DeclareMathOperator{\im}{im}
\DeclareMathOperator{\id}{id}
\newcommand{\FN}{{\ensuremath \F\spann{N}}}

\newcommand{\FcalH}{{\ensuremath \F\spann{\calH}}}
\newcommand{\Fb}{{\ensuremath \F,b}}
\newcommand{\Gb}{{\ensuremath G,b}}
\DeclareMathOperator{\Cone}{Cone}
\DeclareMathOperator{\cd}{cd}
\DeclareMathOperator{\can}{can}

\begin{document}
\maketitle

\begin{abstract}
	We introduce a bounded version of Bredon cohomology for groups relative to a family of subgroups. Our theory generalizes bounded cohomology and differs from Mineyev--Yaman's relative bounded cohomology for pairs. We obtain cohomological characterizations of relative amenability and relative hyperbolicity, analogous to the results of Johnson and Mineyev for bounded cohomology. 
\end{abstract}

\section{Introduction}

Bounded cohomology is a homotopy invariant of topological spaces with deep connections to Riemannian geometry via the simplicial volume of manifolds~\cite{Gromov82}. An astonishing phenomenon known as Gromov's Mapping Theorem is that for any CW-complex $X$, the classifying map $X\to B\pi_1(X)$ induces an isometric isomorphism on bounded cohomology. This emphasizes the importance of the corresponding theory of bounded cohomology for groups, which is also of independent interest due to its plentiful applications in geometric group theory~\cite{Monod01,Monod06,Frigerio17}. 
The bounded cohomology $H^n_b(G;V)$ of a (discrete) group $G$ with coefficients in a normed $G$-module $V$ is the cohomology of the cochain complex of bounded $G$-maps $G^{n+1}\to V$. The inclusion of bounded $G$-maps into (not necessarily bounded) $G$-maps induces the so called comparison map $H^n_b(G;V)\to H^n(G;V)$. On the one hand, the bounded cohomology groups are very difficult to compute in general. On the other hand, they characterize interesting group-theoretic properties such as amenability~\cite{Johnson72} and hyperbolicity~\cite{Mineyev01,Mineyev02}.

\begin{thm}[Johnson]\label{thm:Johnson}
	Let $G$ be a group. The following are equivalent:
	\begin{enumerate}[label=\enum]
		\item $G$ is amenable;
		\item $H^n_b(G;V^\#)=0$ for any dual normed $\IR G$-module $V^\#$ and all $n\ge 1$;
		\item $H^1_b(G;V^\#)=0$ for any dual normed $\IR G$-module $V^\#$.
	\end{enumerate}
\end{thm}

\begin{thm}[Mineyev]\label{thm:Mineyev}
	Let $G$ be a finitely presented group. The following are equivalent:
	\begin{enumerate}[label=\enum]
		\item $G$ is hyperbolic;
		\item The comparison map $H^n_b(G;V)\to H^n(G;V)$ is surjective for any normed $\IQ G$-module $V$ and all $n\ge 2$;
		\item The comparison map $H^2_b(G;V)\to H^2(G;V)$ is surjective for any normed $\IR G$-module $V$.
	\end{enumerate}
\end{thm}

There are well-studied notions of relative amenability and relative hyperbolicity in the literature~\cite{JOR12,Hruska10}. 
In the present article we introduce a new ``relative bounded cohomology theory" characterizing these relative group-theoretic properties as a bounded version of Bredon cohomology. For a group $G$, a family of subgroups $\F$ is a non-empty set of subgroups which is closed under conjugation and taking subgroups. 
For a set of subgroup $\calH$ of $G$, we denote by $\FcalH$ the smallest family containing $\calH$. 
The Bredon cohomology $H^n_\F(G;V)$ with coefficients in a $G$-module $V$ (or more general coefficient systems) is a generalization of group cohomology, which is recovered when $\F$ consists only of the trivial subgroup.
A fundamental feature of Bredon cohomology is that for a normal subgroup $N$ of $G$ there is an isomorphism $H^n_\FN(G;V)\cong H^n(G/N;V^N)$.
From a topological point of view, the Bredon cohomology of $G$ can be identified with the equivariant cohomology of the classifying space $\EFG$ for the family $\F$. Especially the classifying spaces $E_\FIN G$ and $E_\VCY G$ for the family of finite groups and virtually cyclic groups have received a lot of attention in recent years due to their prominent role in the Isomorphism Conjectures of Baum--Connes and Farrell--Jones, respectively.

We introduce the \emph{bounded Bredon cohomology} $H^n_\Fb(G;V)$ of $G$ with coefficients in a normed $G$-module $V$, which generalizes bounded cohomology (Definition~\ref{defn:BBC}). Our theory still is well-behaved with respect to normal subgroups (Corollary~\ref{cor:normal}) and admits a topological interpretation in terms of classifying spaces for families (Theorem~\ref{thm:comb top}). We obtain the following generalizations of Theorems~\ref{thm:Johnson} and~\ref{thm:Mineyev}. A group $G$ is called amenable relative to a set of subgroups $\calH$ if there exists a $G$-invariant mean on the $G$-set $\coprod_{H\in\calH}G/H$.

\begin{thm}\label{thm:thm A}
	Let $G$ be a group and $\calH$ be a set of subgroups. The following are equivalent:
	\begin{enumerate}[label=\enum]
		\item $G$ is amenable relative to $\calH$;
		\item $H^n_{\FcalH,b}(G;V^\#)=0$ for any dual normed $\IR G$-module $V^\#$ and all $n\ge 1$;
		\item $H^1_{\FcalH,b}(G;V^\#)=0$ for any dual normed $\IR G$-module $V^\#$.
	\end{enumerate}
\end{thm}

We also provide a characterization of relative amenability in terms of relatively injective modules (Proposition~\ref{prop:rel amen rel inj}).
Recall that a finite set of subgroups $\calH$ is called a malnormal (resp.\ almost malnormal) collection if for all $H_i,H_j\in\calH$ and $g\in G$ we have $H_i\cap gH_jg^{-1}$ is trivial (resp.\ finite) or $i=j$ and $g\in H_i$. A group $G$ is said to be of type $F_{n,\F}$ for a family of subgroups $\F$, if there exists a model for $\EFG$ with cocompact $n$-skeleton. 

\begin{thm}[Theorem \ref{thm:rel hyp}]\label{thm:thm B}
	Let $G$ be a finitely generated torsionfree group and $\calH$ be a finite malnormal collection of subgroups. Suppose that $G$ is of type $F_{2,\FcalH}$ (e.g.\ $G$ and all subgroups in $\calH$ are finitely presented). Then the following are equivalent:
	\begin{enumerate}[label=\enum]
		\item\label{item:thm B i} $G$ is hyperbolic relative to $\calH$;
		\item The comparison map $H^n_{\FcalH,b}(G;V)\to H^n_\FcalH(G;V)$ is surjective for any normed $\IQ G$-module $V$ and all $n\ge 2$;
		\item\label{item:thm B iii} The comparison map $H^2_{\FcalH,b}(G;V)\to H^2_\FcalH(G;V)$ is surjective for any normed $\IR G$-module $V$.
	\end{enumerate}
\end{thm}
In Theorem~\ref{thm:thm B} the equivalence of~\ref{item:thm B i} and~\ref{item:thm B iii} still holds if the group $G$ contains torsion and $\calH$ is almost malnormal, see Remark~\ref{rem:torsion}.
Note that condition~\ref{item:thm B iii} is trivially satisfied for groups of Bredon cohomological dimension $\cd_{\FcalH}$ equal to 1.

The topological interpretation of bounded Bredon cohomology via classifying spaces for families was used by L\"oh--Sauer~\cite{Loeh-Sauer19} to give a new proof of the Nerve Theorem and Vanishing Theorem for amenable covers. 
We prove a converse of~\cite[Proposition 5.2]{Loeh-Sauer19}, generalizing a recent result of~\cite[Theorem 3.1.3]{Moraschini-Raptis21} where the case of a normal subgroup is treated.

\begin{thm}\label{thm:thm C}
	Let $G$ be a group and $\F$ be a family of subgroups. The following are equivalent:
	\begin{enumerate}[label=\enum]
		\item All subgroups in $\F$ are amenable;
		\item The canonical map $H^n_\Fb(G;V^\#)\to H^n_b(G;V^\#)$ is an isomorphism for any dual normed $\IR G$-module $V^\#$ and all $n\ge 0$;
		\item The canonical map $H^1_\Fb(G;V^\#)\to H^1_b(G;V^\#)$ is an isomorphism for any dual normed $\IR G$-module $V^\#$.
	\end{enumerate}
\end{thm}
In fact, both Theorems~\ref{thm:thm A} and~\ref{thm:thm C} are special cases of the more general Theorem~\ref{thm:unification}.
As an application of Theorem~\ref{thm:thm C}, the comparison map vanishes for groups which admit a ``small" model for $\EFG$, where $\F$ is any family consisting of amenable subgroups (Corollary~\ref{cor:Loeh-Sauer}). Examples are graph products of amenable groups (e.g.\ right-angled Artin groups) and fundamental groups of graphs of amenable groups.

There is another natural relative cohomology theory given by the relative cohomology of a pair of spaces. For a set of subgroups $\calH$, it gives rise to the cohomology $H^n(G,\calH;V)$ of the group pair $(G,\calH)$ introduced by Bieri--Eckmann~\cite{Bieri-Eckmann78}. A bounded version $H^n_b(G,\calH;V)$ was defined by Mineyev--Yaman~\cite{Mineyev-Yaman07} to give a characterization of relative hyperbolicity (see also~\cite{Franceschini18}). A characterization of relative amenability in terms of this relative theory was obtained in~\cite{JOR12}. There is a canonical map $H^n_\FcalH(G;V)\to H^n(G,\calH;V)$ for $n\ge 2$ which is an isomorphism if $\calH$ is malnormal (see Remark~\ref{rem:Bieri-Eckmann}). Similarly, there is a map for the bounded versions but we do not know when it is an isomorphism due to the failure of the excision axiom for bounded cohomology (see Remark~\ref{rem:Mineyev-Yaman}). We also mention that Mineyev--Yaman's relative bounded cohomology was extended to pairs of groupoids in~\cite{Blank16}.
\\
\paragraph{Acknowledgements}
The present work is part of the author's PhD project. He wishes to thank his supervisors Nansen Petrosyan for suggesting this topic as well as for numerous discussions and Ian Leary for his support. We are grateful to the organizers of the ``Virtual workshop: Simplicial Volumes and Bounded Cohomology" held in September 2020 during which parts of this work were discussed. We thank Clara L\"oh for several interesting conversations, Francesco Fournier-Facio, Sam Hughes, and Eduardo Mart\'inez-Pedroza for helpful comments.

\tableofcontents

\section{Preliminaries on Bredon cohomology and classifying spaces}

In this section we briefly recall the notion of Bredon cohomology for groups and its topological interpretation as the equivariant cohomology of classifying spaces for families of subgroups. For an introduction to Bredon cohomology we refer to~\cite{Fluch10} and for a survey on classifying spaces to~\cite{Lueck05survey}.

Let $G$ be a group, which shall always mean a discrete group.
A \emph{family of subgroups $\F$} is a non-empty set of subgroups of $G$ that is closed under conjugation by elements of $G$ and under taking subgroups. Typical examples are $\TR=\{1\}$, $\FIN=\{\text{finite subgroups}\}$, $\VCY=\{\text{virtually cyclic subgroups}\}$, and $\ALL=\{\text{all subgroups}\}$. We will moreover be interested in $\AME=\{\text{amenable subgroups}\}$. 
For a subgroup $H$ of $G$, we denote by $\F|_H$ the family $\{L\cap H\ |\ L\in\F\}$ of subgroups of $H$. (In the literature this family is sometimes denoted by $\F\cap H$ instead.)
For a set of subgroups $\calH$, one can consider the smallest family containing $\calH$ which is $\F\spann{\calH}=\{\text{conjugates of elements in $\calH$ and their subgroups}\}$ and called the \emph{family generated by $\calH$}. When $\calH$ consists of a single subgroup $H$, we denote $\F\spann{\calH}$ instead by $\F\spann{H}$ and call it the \emph{family generated by $H$}. We denote by $G/\calH$ the $G$-set $\coprod_{H\in \calH}G/H$.

Let $R$ be a ring and $\mathbf{Mod}_R$ denote the category of $R$-modules. We will often suppress the ring $R$, so that $G$-modules are understood to be $RG$-modules.
The \emph{($\F$-restricted) orbit category $\OFG$} has as objects $G$-sets of the form $G/H$ with $H\in \F$ and as morphisms $G$-maps. 
An \emph{$\OFG$-module} is a contravariant functor $M\colon \OFG\to \mathbf{Mod}_R$, the category of which is denoted by $\OFG\Mod_R$.
Note that $\OTRG\Mod_R$ can be identified with the category of $G$-modules. 
For a $G$-module $V$, there is a coinduced $\OFG$-module $V^?$ given by $V^?(G/H)=V^H$. (In the literature this is sometimes called a \emph{fixed-point functor}.) Observe that $(-)^?$ is right-adjoint to the restriction $\OFG\Mod_R\to \OTRG\Mod_R$, $M\mapsto M(G/1)$.
For a $G$-space $X$ and a $G$-CW-complex $Y$ with stabilizers in $\F$, there are singular and cellular $\OFG$-chain complexes $C_*(X^?)(G/H)=C_*(X^H)$ and $C_*^\text{cell}(Y^?)(G/H)=C_*^\text{cell}(Y^H)$, where $C_*(X^H)$ and $C_*^\text{cell}(Y^H)$ denote the usual singular and cellular chain complexes, respectively.

The \emph{Bredon cohomology} of $G$ with coefficients in an $\OFG$-module $M$ is defined as
\[
	H^n_\F(G;M)\coloneqq \Ext^n_{\OFG\Mod_R}(R,M)
\]
for $n\ge 0$, where $R$ is regarded as a constant $\OFG$-module. It can be computed as the cohomology of the cochain complex $\Hom_{\OFG\Mod_R}(R[((G/\F)^{*+1})^?],M)$.
We define the $G$-chain complex $C_*^\F(G)$ given by $G$-modules
\[
	C_n^\F(G)\coloneqq R[(G/\F)^{n+1}]
\]
with the diagonal $G$-action and differentials $\partial_n\colon C_n^\F(G)\to C_{n-1}^\F(G)$,
\[
	\partial_n(g_0H_0,\ldots,g_nH_n)=\sum_{i=0}^n (-1)^i (g_0H_0,\ldots,\widehat{g_iH_i},\ldots,g_nH_n) \,.
\]
For a $G$-module $V$, the $G$-cochain complex $C^*_\F(G;V)$ is given by
	\[
		C^n_\F(G;V)\coloneqq \Hom_R(C_n^\F(G),V)
	\]
	so that
	\[
		H^n_\F(G;V)\coloneqq H^n_\F(G;V^?)\cong H^n(C^*_\F(G;V)^G) \,.
	\]
For a $G$-space $X$ with stabilizers in $\F$, the \emph{Bredon cohomology} of $X$ with coefficients in an $\OFG$-module $M$ is defined as
\[
	H^n_G(X;M)\coloneqq H^n(\Hom_{\OFG\Mod_R}(C_*(X^?),M))
\]
for $n\ge 0$.
If $X$ is a $G$-CW-complex, then $H^n_G(X;M)$ can be computed using $C_*^\cell(X^?)$ instead of $C_*(X^?)$.

A \emph{classifying space $\EFG$ for the family $\F$} is a terminal object in the $G$-homotopy category of $G$-CW-complexes with stabilizers in $\F$.
It can be shown that a $G$-CW-complex $X$ is a model for $\EFG$ if and only if the fixed-point set $X^H$ is contractible for $H\in \F$ and empty otherwise.
An explicit model is given by the geometric realization $Y$ of the semi-simplicial set $\{(G/\F)^{n+1}\ |\ n\ge 0\}$ with the usual face maps. Then $Y$ has (non-equivariant) $n$-cells corresponding to $(G/\F)^{n+1}$ and we refer to $Y$ as the \emph{simplicial model} for $\EFG$.
Note that a model for $E_\TR G$ is given by $EG$ and a model for $E_\ALL G$ is the point $G/G$.
The cellular $\OFG$-chain complex of any model for $\EFG$ is a projective resolution of the constant $\OFG$-module $R$ and thus we have
\begin{equation}\label{eqn:Bredon EFG}
	H^n_\F(G;M)\cong H^n_G(\EFG;M)
\end{equation}
for any $\OFG$-module $M$.
If $N$ is a normal subgroup of $G$, then a model for $E_{\F\spann{N}}G$ is given by $E(G/N)$ regarded as a $G$-CW-complex and we find
\begin{equation}\label{eqn:Bredon quotient}
	H^n_{\F\spann{N}}(G;M)\cong H^n(G/N;M(G/N)) \,.
\end{equation}

For a subgroup $H$ of $G$, when viewed as an $H$-space $\EFG$ is a model for $E_{\F|_H}H$ which induces the restriction map
\begin{equation}\label{eqn:Bredon res}
	\res^n_{H\subset G}\colon H^n_\F(G;M)\to H^n_{\F|_H}(H;M)
\end{equation}
for any $\OFG$-module $M$.
For two families of subgroups $\F\subset \G$, the up to $G$-homotopy unique $G$-map $\EFG\to \EGG$ induces the canonical map
\begin{equation}\label{eqn:Bredon can}
	\can^n_{\F\subset\G}\colon H^n_\G(G;M)\to H^n_\F(G;M)
\end{equation}
for any $\mathcal{O}_\G G$-module $M$.

\begin{rem}[Bieri--Eckmann's relative cohomology]\label{rem:Bieri-Eckmann}
	For a group $G$ and a set of subgroups $\calH$, Bieri--Eckmann~\cite{Bieri-Eckmann78} have introduced the \emph{relative cohomology $H^n(G,\calH;V)$ of the pair $(G,\calH)$} with coefficients in a $G$-module $V$. It can be identified with the relative cohomology $H^n_G(EG,\coprod_{H\in \calH}G\times_H EH;V)$ of the pair of $G$-spaces $(EG,\coprod_{H\in\calH}G\times_H EH)$. Here a model for $EG$ is chosen that contains $\coprod_{H\in\calH} G\times_H EH$ as a subcomplex by taking mapping cylinders.
	Hence there is a long exact sequence
	\[
		\cdots H^n(G,\calH;V)\to H^n(G;V)\to \prod_{H\in\calH}H^n(H;V)\to \cdots \,,
	\]
	which is one of the main features of the relative cohomology groups.
	
	There is a relation between Bredon cohomology and Bieri--Eckmann's relative cohomology as follows.
	Consider the $G$-space $X$ obtained as the $G$-pushout
	\[\begin{tikzcd}
		\coprod_{H\in\calH} G\times_H EH\ar{r}\ar{d} & EG\ar{d} \\
		\coprod_{H\in\calH} G/H\ar{r} & X \,,
	\end{tikzcd}\]
	where the left vertical map is induced by collapsing each $EH$ to a point.
	Then the $G$-space $X$ has stabilizers in $\FcalH$ and hence admits a $G$-map $X\to E_\FcalH G$. For an $\OFG$-module $M$, we have maps
	\[\begin{tikzcd}
		H^n_G(X;M) & H^n_G(X,\coprod_{H\in\calH}G/H;M)\ar{l}\ar{d}{\cong} \\
		H^n_G(E_\FcalH G;M)\ar{u} & H^n_G(EG,\coprod_{H\in\calH}G\times_H EH;M) \,,
	\end{tikzcd}\]
	where the right vertical map is an isomorphism by excision.
	Now, if $\calH$ is a malnormal collection, then $X$ is a model for $E_{\F\spann{\calH}}G$ and we have	
	\[
		H^n_{\F\spann{\calH}}(G;M)\cong H^n(G,\calH;M(G/1))
	\]
	for $n\ge 2$. This was shown in~\cite[Theorem 4.16]{AN-CM-SS18} for the special case when $\calH$ consists of a single subgroup.
\end{rem}

\section{Bounded Bredon cohomology}
In this section we introduce a bounded version of Bredon cohomology and develop some of its basic properties. We follow the exposition in~\cite{Frigerio17} for bounded cohomology. Throughout, let $G$ be a group and $\F$ be a family of subgroups.

From now on, let the ring $R$ be one of $\IZ$, $\IQ$ or $\IR$. A \emph{normed $G$-module} $V$ is a $G$-module equipped with a $G$-invariant norm $\|\cdot\|\colon V\to \IR$. (That is, for all $v,u\in V$, $r\in R$, and $g\in G$ we have $\|v\|=0$ if and only if $v=0$, $\|rv\|\le |r|\cdot \|v\|$, $\|v+u\|\le \|v\|+\|u\|$, and $\|g\cdot v\|=\|v\|$.) A morphism $f\colon V\to W$ of normed $G$-modules is a morphism of $G$-modules with finite operator-norm $\|f\|_\infty$. By $\bHom_R(V,W)$ we denote the $G$-module of $R$-linear maps $f\colon V\to W$ with finite operator-norm, where the $G$-action is given by $(g\cdot f)(v)=g\cdot f(g^{-1}v)$. We denote the topological dual $\bHom_R(V,\IR)$ of $V$ by $V^\#$. For a set $S$ and a normed module $V$, we denote by $\bMap(S,V)$ the module of functions $S\to V$ with bounded image. Instead of $\bMap(S,\IR)$ we also write $\linf(S)$.

The following is our key definition.
Recall the notation $G/\F = \coprod_{H\in \F}G/H$ and consider $C_n^\F(G)=R[(G/\F)^{n+1}]$ as a normed $G$-module equipped with the $\ell^1$-norm with respect to the generators. For a normed $G$-module $V$, we define the cochain complex $C^*_{\F,b}(G;V)$ of normed $G$-modules by
\[
	C^n_{\F,b}(G;V)\coloneqq \bHom_R(C_n^\F(G),V)
\]
together with the differentials $\delta^n\colon C^n_\Fb(G;V)\to C^{n+1}_\Fb(G;V)$,
\[
	\delta^n(f)(g_0H_0,\ldots,g_{n+1}H_{n+1})=\sum_{i=0}^{n+1}(-1)^i f(g_0H_0,\ldots,\widehat{g_iHi},\ldots,g_{n+1}H_{n+1}) \,.
\]

\begin{defn}[Bounded Bredon cohomology of groups]\label{defn:BBC}
The \emph{bounded Bredon cohomology} of $G$ with coefficients in a normed $G$-module $V$ is defined as
\[
	H^n_{\F,b}(G;V)\coloneqq H^n(C^*_{\F,b}(G;V)^G)
\]
for $n\ge 0$.
The inclusion $C^n_\Fb(G;V)\subset C^n_\F(G;V)$ induces a map
\[
	c^n_\F\colon H^n_\Fb(G;V)\to H^n_\F(G;V)
\]
called the \emph{comparison map}.
\end{defn}

Note that for $\F=\TR$, Definition~\ref{defn:BBC} recovers the usual definition of bounded cohomology.
\begin{rem}[Coefficient modules]
	We only consider normed $G$-modules as coefficients, rather than more general $\OFG$-modules equipped with a ``compatible norm". 
	Hence strictly speaking our theory is a bounded version of Nucinkis' cohomology relative to the $G$-set $G/\F$~\cite{Nucinkis99}, rather than a bounded version of Bredon cohomology.
\end{rem}

\begin{rem}[Canonical semi-norm]
	The $\ell^\infty$-norm on $C^n_{\F,b}(G;V)$ descends to a \emph{canonical semi-norm} on $H^n_{\F,b}(G;V)$. However, we do not consider semi-norms anywhere in this article and regard $H^n_{\F,b}(G;V)$ merely as an $R$-module. 
\end{rem}

Bounded Bredon cohomology satisfies the following basic properties.
\begin{lem}\label{lem:basic properties}
	The following hold:
	\begin{enumerate}[label=\enum]
		\item\label{item:basic i} Let $0\to V_0\to V_1\to V_2\to 0$ be a short exact sequence of normed $G$-modules such that $0\to V_0^H\to V_1^H\to V_2^H\to 0$ is exact for each $H\in\F$. Then there exists a long exact sequence
		\[
			0\to H^0_{\F,b}(G;V_0)\to H^0_{\F,b}(G;V_1)\to H^0_{\F,b}(G;V_2)\to H^1_{\F,b}(G;V_0)\to \cdots \,;
		\] 
		\item\label{item:basic ii} $H^0_\Fb(G;V)\cong V^G$ for any normed $G$-module $V$;
		\item $H^1_\Fb(G;\IR)=0$.
	\end{enumerate}
	\begin{proof}
		(i) Observe that the sequence of cochain complexes
		\[
			0\to C^*_\Fb(G;V_0)^G\to C^*_\Fb(G;V_1)^G\to C^*_\Fb(G;V_2)^G\to 0
		\]
		is exact. Then the associated long exact sequence on cohomology is as desired.
		
		(ii) We have $H^0_\Fb(G;V)=\ker(\delta^0)$, where 
		\[
			\delta^0\colon \bHom_{RG}(R[G/\F],V)\to \bHom_{RG}(R[(G/\F)^2],V)
		\]
		is given by $\delta^0(f)(g_0H_0,g_1H_1)=f(g_1H_1)-f(g_0H_0)$. Hence $\ker(\delta^0)$ consists precisely of the constant $G$-maps $G/\F\to V$, which are in correspondence to $V^G$.
		
		(iii) We identify 
		\[
			C^n_\Fb(G;\IR)^G\cong \bMap\big(\coprod_{H_0,\ldots,H_n\in\F}H_0\backslash(G/H_1\times\cdots \times G/H_n),\IR\big)
		\]
		for $n\ge 1$ and $C^0_\Fb(G;\IR)^G\cong \bMap\left(\coprod_{H_0\in\F}\ast_{H_0},\IR\right)$. The differentials of this ``inhomogeneous" complex in low degrees are given by
		\begin{align*}
			\delta^0(f)(H_0g_1H_1) &= f(\ast_{H_1})-f(\ast_{H_0}) \\
			\delta^1(\varphi)(H_0(g_1H_1,g_2H_2)) &= \varphi(H_1g_1^{-1}g_2H_2)-\varphi(H_0g_2H_2)+\varphi(H_0g_1H_1) \,.
		\end{align*}
		Then it is not difficult to check that $\ker(\delta^1)=\im(\delta^0)$.
	\end{proof}
\end{lem}

	We also define the bounded cohomology of a $G$-space $X$ as follows. Denote by $S_n(X)$ the set of singular $n$-simplices in $X$ and consider $C_n(X)=R[S_n(X)]$ equipped with the $\ell^1$-norm as a normed $G$-module. For a normed $G$-module $V$, we define the cochain complex $C^*_b(X;V)$ of normed $G$-modules by
	\[
		C^n_b(X;V)\coloneqq \bHom_R(C_n(X),V)
	\]
	together with the usual differentials. 
	
\begin{defn}[Bounded cohomology of $G$-spaces]	
	The \emph{($G$-equivariant) bounded cohomology} of a $G$-space $X$ with coefficients in a normed $G$-module $V$ is defined as
	\[
		H^n_\Gb(X;V)\coloneqq H^n(C^*_b(X;V)^G)
	\]
	for $n\ge 0$. The inclusion $C^n_b(X;V)\subset C^n(X;V)$ induces a map
	\[
		c^n_X\colon H^n_\Gb(X;V)\to H^n_G(X;V)
	\]
	called the \emph{comparison map}.
\end{defn}
Note that the functors $H^*_\Gb$ are $G$-homotopy invariant and that $H^n_\Gb(G/H;V)$ is isomorphic to $V^H$ for $n=0$ and trivial otherwise. However, beware that $H^*_\Gb$ is neither a $G$-cohomology theory, nor can it be computed cellularly for $G$-CW-complexes, as is the case already when $G$ is the trivial group. 
\\
\paragraph{Relative homological algebra}
We develop the relative homological algebra that will allow us to compute bounded Bredon cohomology via resolutions, analogous to Ivanov's approach for bounded cohomology~\cite{Ivanov87}.

	A map $p\colon A\to B$ of $G$-modules is called \emph{$\F$-strongly surjective} if for each $H\in\F$ there exists a map $\tau_H\colon B\to A$ of $H$-modules such that $p\circ \tau_H=\id_B$. A $G$-module $P$ is called \emph{relatively $\F$-projective} if for every $\F$-strongly surjective $G$-map $p\colon A\to B$ and every $G$-map $\phi\colon P\to B$, there exists a $G$-map $\Phi\colon P\to A$ such that $p\circ \Phi=\phi$.
	A chain complex of $G$-modules is called \emph{relatively $\F$-projective} if each chain module is relatively $\F$-projective.
	A resolution $(C_*,\partial_*)$ of $G$-modules is called \emph{$\F$-strong} if it is contractible as a resolution of $H$-modules for each $H\in\F$. (That is, there exist $H$-maps $k^H_*\colon C_*\to C_{*+1}$ such that $\partial_{n+1}\circ k^H_n+k^H_{n-1}\circ \partial_n=\id_{C_n}$.)
	
\begin{lem}\label{lem:rel proj}
	The following hold:
	\begin{enumerate}[label=\enum]
		\item If $S$ is a $G$-set with stabilizers in $\F$, then the $G$-module $R[S]$ is relatively $\F$-projective;
		\item If $S$ is a $G$-set with $S^H\neq \emptyset$ for all $H\in\F$, then the resolution $R[S^{*+1}]\to R$ of $G$-modules is $\F$-strong;
		\item\label{item:rel proj G-space} If $X$ is a $G$-space with contractible fixed-point set $X^H$ for each $H\in \F$, then the resolution $C_*(X)\to R$ of $G$-modules is $\F$-strong.
	\end{enumerate}
	\begin{proof}
		(i) Given a lifting problem as in the definition of relative $\F$-projectivity, 
		\[\begin{tikzcd}
			& R[S]\ar{d}{\phi}\ar[dashed]{dl}[swap]{\Phi} \\
			A\ar{r}{p} & B\ar{r}\ar[bend left]{l}{\tau_H} & 0
		\end{tikzcd}\]
		we construct a lift $\Phi$ as follows. Let $T$ be a set of representatives of $G\backslash S$ and denote the stabilizer of an element $t\in T$ by $G_t$. Then for every $s\in S$ there exist unique elements $t_s\in T$ and $g_sG_{t_s}\in G/G_{t_s}$ such that $g_s^{-1}s=t_s$. Define $\Phi\colon R[S]\to A$ on generators by
		\[
			\Phi(s)=g_s\cdot \tau_{G_{t_s}}(\phi(g_s^{-1}s))
		\]
		which is a well-defined $G$-equivariant lift of $\phi$.
		
		(ii) For $H\in\F$, fix an element $s_H\in S^H$ and define $k^H_*\colon R[S^{*+1}]\to R[S^{*+2}]$ on generators by
		\[
			k^H_n(s_0,\ldots,s_n)=(s_H,s_0,\ldots,s_n) \,.
		\]
		Then $k^H_*$ is an $H$-equivariant contraction.
		
		(iii) For $H\in\F$, fix a point $x_H\in X^H$ and define a contraction $k_*^H\colon C_*(X)\to C_{*+1}(X)$ of $H$-chain complexes inductively as follows. Starting with $k_{-1}^H\colon R\to C_0(X)$, $r\mapsto r\cdot x_H$, we may assume that $k_{n-1}^H$ has been constructed. Let $s$ be a singular $n$-simplex in $X$ and denote its stabilizer by $H_s$. Then there exists a singular $(n+1)$-simplex $s'$ with 0-th vertex $x_H$ and opposite face $s$, satisfying $\partial_{n+1}(s')+k_{n-1}^H(\partial_n(s))=s$. Moreover, since $X^{H_s}$ is contractible we may choose $s'$ such that its image is contained in $X^{H_s}$. Now, for each $H$-orbit of singular $n$-simplices in $X$ choose a representative $s$, define $k_n^H(s)$ to be $s'$ and then extend $H$-equivariantly.
	\end{proof}
\end{lem}	

The proof of the following proposition is standard and omitted. 

\begin{prop}\label{prop:fund lemma I}
	Let $f\colon V\to W$ be a map of $G$-modules, $P_*\to V$ be a $G$-chain complex with $P_n$ relatively $\F$-projective for all $n\ge 0$, and $C_*\to W$ be an $\F$-strong resolution of $G$-modules. Then there exists a $G$-chain map $f_*\colon P_*\to C_*$ extending $f$, which is unique up to $G$-chain homotopy.
\end{prop}

	While relatively $\F$-projective $\F$-strong resolutions are useful to compute Bredon homology, the following dual approach will compute bounded Bredon cohomology.

	A map $i\colon A\to B$ of normed $G$-modules is called \emph{$\F$-strongly injective} if for each $H\in\F$ there exists a map $\sigma_H\colon B\to A$ of normed $H$-modules with $\|\sigma_H\|_\infty\le K$ such that $\sigma_H\circ i=\id_A$, for a uniform constant $K\ge 0$. 
	A normed $G$-module $I$ is called \emph{relatively $\F$-injective} if for every $\F$-strongly injective $G$-map $i\colon A\to B$ and every map $\psi\colon A\to I$ of normed $G$-modules, there exists a map $\Psi\colon B\to I$ of normed $G$-modules such that $\Psi\circ i=\psi$.
	A chain complex of normed $G$-modules is called \emph{relatively $\F$-injective} if each chain module is relatively $\F$-injective.
	A resolution of normed $G$-modules is called \emph{$\F$-strong} if it is contractible as a resolution of normed $H$-modules for each $H\in\F$.

Dual to Lemma~\ref{lem:rel proj} and Proposition~\ref{prop:fund lemma I} we obtain the following.
\begin{lem}\label{lem:ex rel inj}
	Let $V$ be a normed $G$-module. The following hold:
	\begin{enumerate}[label=\enum]
		\item If $S$ is a $G$-set with stabilizers in $\F$, then $\bHom_R(R[S],V)$ is a relatively $\F$-injective normed $G$-module;
		\item\label{item:F-strong} If $S$ is a $G$-set with $S^H\neq \emptyset$ for all $H\in\F$, then the resolution $V\to \bHom_R(R[S^{*+1}],V)$ of normed $G$-modules is $\F$-strong;
		\item If $X$ is a $G$-space with contractible fixed-point set $X^H$ for each $H\in \F$, then the resolution $V\to C^*_b(X;V)$ of normed $G$-modules is $\F$-strong.
	\end{enumerate}
\end{lem}

\begin{prop}\label{prop:fund lemma II}
	Let $f\colon V\to W$ be a map of normed $G$-modules, $V\to C^*$ be an $\F$-strong resolution of normed $G$-modules, and $W\to I^*$ be a $G$-chain complex with $I^n$ relatively $\F$-injective for all $n\ge 0$. Then there exists a $G$-chain map $f^*\colon C^*\to I^*$ extending $f$, which is unique up to $G$-chain homotopy.
\end{prop}

As a consequence of Proposition~\ref{prop:fund lemma II}, we may use any relatively $\F$-injective $\F$-strong resolution to compute bounded Bredon cohomology. We obtain the isomorphisms analogous to~\eqref{eqn:Bredon EFG} and~\eqref{eqn:Bredon quotient} for Bredon cohomology.

\begin{thm}\label{thm:comb top}
	Let $G$ be a group, $\F$ be a family of subgroups, and $V$ be a normed $G$-module. For all $n\ge 0$ there is an isomorphism
	\[
		H^n_\Fb(G;V)\cong H^n_\Gb(\EFG;V) \,.
	\]
	\begin{proof}
		Both $C^*_\Fb(G;V)$ and $C^*_b(\EFG;V)$ are relatively $\F$-injective $\F$-strong resolutions of $V$ by Lemma~\ref{lem:ex rel inj} and hence $G$-chain homotopy equivalent by Proposition~\ref{prop:fund lemma II}.
	\end{proof}
\end{thm}

\begin{cor}\label{cor:normal}
	Let $G$ be a group, $N$ be a normal subgroup of $G$, and $V$ be a normed $G$-module. For all $n\ge 0$ there is an isomorphism
	\[
		H^n_{\FN,b}(G;V)\cong H^n_b(G/N;V^N) \,.
	\]
	\begin{proof}
		As a model for $E_\FN G$ we take $E(G/N)$ regarded as a $G$-space. Then it suffices to observe that
		\[
			\bHom_{RG}(R[S_n(E(G/N))],V) \cong \bHom_{R[G/N]}(R[S_n(E(G/N))],V^N)
		\]
		and to apply Theorem~\ref{thm:comb top} twice.
	\end{proof}
\end{cor}

	Analogous to~\eqref{eqn:Bredon res} and~\eqref{eqn:Bredon can} for Bredon cohomology, for a subgroup $H$ of $G$ and two families of subgroups $\F\subset \G$, we have the maps
	\begin{align*}
		\res^n_{H\subset G,b}\colon & H^n_\Fb(G;V)\to H^n_{\F|_H,b}(H;V) \,; \\
		\can^n_{\F\subset \G,b}\colon & H^n_{\G,b}(G;V)\to H^n_\Fb(G;V)
	\end{align*}
	for any normed $G$-module $V$.

\begin{rem}[Mineyev--Yaman's relative bounded cohomology]\label{rem:Mineyev-Yaman}
	Mineyev--Yaman have introduced the bounded analogue of Bieri--Eckmann's relative cohomology for pairs (Remark~\ref{rem:Bieri-Eckmann}) in~\cite{Mineyev-Yaman07}. For a group $G$, a set of subgroups $\calH$, and a normed $G$-module $V$, their relative bounded cohomology groups $H^n_b(G,\calH;V)$ can be identified with $H^n_\Gb(EG,\coprod_{H\in\calH}G\times_H EH;V)$ and therefore fit in a long exact sequence
	\[
		\cdots\to H^n_b(G,\calH;V)\to H^n_b(G;V)\to \prod_{H\in\calH} H^n_b(H;V)\to \cdots \,.
	\]
	As in Remark~\ref{rem:Bieri-Eckmann}, we denote by $X$ the $G$-space obtained as a $G$-pushout from $EG$ by collapsing $G\times_H EH$ to $G/H$ for each $H\in\calH$. Then we have maps
	\[\begin{tikzcd}
		H^n_\Gb(X;V) & H^n_\Gb(X,\coprod_{H\in\calH}G/H;V)\ar{l}\ar{d} \\		
		H^n_\Gb(E_\FcalH G;V)\ar{u} & H^n_\Gb(EG,\coprod_{H\in\calH}G\times_H EH;V) \,.
	\end{tikzcd}\]
	For $n\ge 2$, the horizontal map is an isomorphism and hence we obtain a map
	\[
		H^n_{\FcalH,b}(G;V)\to H^n_b(G,\calH;V) \,.
	\]
	However, even if $\calH$ is a malnormal collection in which case $X$ is a model for $E_\FcalH G$, this map need not be an isomorphism due to the failure of the excision axiom for bounded cohomology.
\end{rem}

\section{Characterization of relative amenability}\label{sec:rel amen}

In this section we prove a characterization of relatively amenable groups in terms of bounded Bredon cohomology analogous to Theorem~\ref{thm:Johnson}.

Recall that a \emph{$G$-invariant mean} on a $G$-set $S$ is an $\IR$-linear map $m\colon \linf(S)\to \IR$ which is normalized, non-negative, and $G$-invariant. (That is, for the constant function $1\in \linf(S)$, $f\in \linf(S)$, and $g\in G$ we have $m(1)=1$, $m(f)\ge 0$ if $f\ge 0$, and $m(g\cdot f)=m(f)$.) 
Note that for a $G$-map $S_1\to S_2$ of $G$-sets, a $G$-invariant mean on $S_1$ is pushed forward to a $G$-invariant mean on $S_2$.

\begin{defn}[Relative amenability]
	 A group $G$ is \emph{amenable relative} to a set of subgroups $\calH$ if the $G$-set $G/\calH$ admits a $G$-invariant mean. When $G$ is amenable relative to $\calH$ consisting of a single subgroup $H$, we also say that $H$ is \emph{co-amenable} in $G$. 
\end{defn}

When $\calH$ is a finite set of subgroups, we recover the notion of relative amenability studied in~\cite{JOR12} (see also~\cite{Monod-Popa03}). 

\begin{ex}
	Let $G$ be a group, $H$ be a subgroup, and $\calH$ be a set of subgroups.	
	\begin{enumerate}[label=\enum]
		\item If $G$ is amenable, then $G$ is amenable relative to $\calH$;
		\item If $H$ is a normal subgroup, then $H$ is co-amenable in $G$ if and only if the quotient group $G/H$ is amenable;
		\item If $H$ has finite index in $G$ or contains the commutator subgroup $[G,G]$, then $H$ is co-amenable in $G$;
		\item If $\calH$ is finite and $G$ is amenable relative to $\calH$, then $\calH$ contains an element that is co-amenable in $G$;
		\item $G$ is amenable relative to $\calH$ if and only if $G$ is amenable relative to $\FcalH$.
	\end{enumerate}
\end{ex}

The following lemma is proved analogously to~\cite[Lemma 3.2]{Frigerio17}.
\begin{lem}\label{lem:rel amen equiv}
	Let $G$ be a group and $\calH$ be a set of subgroups. 
	Then $G$ is amenable relative to $\calH$ if and only if there exists a non-trivial $G$-invariant element in $\linf(G/\calH)^\#$.
\end{lem}

By Proposition~\ref{prop:fund lemma II} bounded Bredon cohomology can be computed using relatively $\F$-injective $\F$-strong resolutions. If one considers coefficients in dual normed $\IR G$-modules, then such resolutions can be obtained from $G$-sets whose stabilizers are amenable relative to $\F$.

\begin{lem}\label{lem:amen rel inj}
	Let $G$ be a group, $\F$ be a family of subgroups, and $V^\#$ be a dual normed $\IR G$-module. If $S$ is a $G$-set such that any stabilizer $G_s$ is amenable relative to $\F|_{G_s}$, then the normed $\IR G$-module $\bHom_R(R[S^{n+1}],V^\#)$ is relatively $\F$-injective for all $n\ge 0$.
	\begin{proof}
		Since the stabilizers of $S^{n+1}$ are intersections of stabilizers of $S$, and relative amenability passes to subgroups, it is enough to consider the case $n=0$.
		Let an extension problem as in the definition of relative $\F$-injectivity be given. 
		\[\begin{tikzcd}
			0\ar{r} & A\ar{r}[swap]{i}\ar{d}[swap]{\psi} & B\ar[dashed]{dl}{\Psi}\ar[bend right]{l}[swap]{\sigma_H} \\
			& \bHom_R(R[S],V^\#)
		\end{tikzcd}\]
		Let $T$ be a set of representatives of $G\backslash S$. We denote the stabilizer of an element $t\in T$ by $G_t$ and by assumption there exists a $G_t$-invariant mean $m_t$ on $G_t/\F|_{G_t}$. 
		Note that any subgroup $L\in \F|_{G_t}$ can also be viewed as an element in $\F$.
		Now, for every $s\in S$ there exist unique elements $t_s\in T$ and $g_sG_{t_s}\in G/G_{t_s}$ such that $g_s^{-1}s=t_s$. Define $\Psi\colon B\to \bHom_R(R[S],V^\#)$ for $b\in B$, $s\in S$, and $v\in V$ by
		\[
			\Psi(b)(s)(v)=m_{t_s}\big(gL\mapsto (g_sg\cdot \psi(\sigma_L(g^{-1}g_s^{-1}b)))(s)(v)\big) \,.
		\]
		One checks that $\Psi$ is a well-defined map of normed $\IR G$-modules extending $\psi$.
	\end{proof}
\end{lem}

For a family of subgroups $\F$, consider the short exact sequence of normed $\IR G$-modules
\[
	0\to \IR\to \linf(G/\F)\to \linf(G/\F)/\IR\to 0 \,,
\]
where $\IR$ is regarded as the constant functions, and the exact sequence of topological duals
\[
	0\to (\linf(G/\F)/\IR)^\#\to \linf(G/\F)^\#\to \IR\to 0 \,.
\]
We define the \emph{relative Johnson class} $[J_\F]\in H^1_\Fb(G;(\linf(G/\F)/\IR)^\#)$ as the cohomology class of the 1-cocycle $J_\F\in C^1_\Fb(G;(\linf(G/\F)/\IR)^\#)$ given by 
\[
	J_\F(g_0H_0,g_1H_1)=\epsilon_{g_1H_1}-\epsilon_{g_0H_0} \,,
\]
where $\epsilon_{g_iH_i}$ is the evaluation map at $g_iH_i$ for $i=0,1$. 

\begin{thm}\label{thm:unification}
	Let $G$ be a group and $\F\subset \G$ be two families of subgroups. The following are equivalent:
	\begin{enumerate}[label=\enum]
		\item Any subgroup $H\in\G$ is amenable relative to $\F|_H$;
		\item The canonical map $H^n_{\G,b}(G;V^\#)\to H^n_\Fb(G;V^\#)$ is an isomorphism for any dual normed $\IR G$-module $V^\#$ and all $n\ge 0$;
		\item The canonical map $H^1_{\G,b}(G;V^\#)\to H^1_\Fb(G;V^\#)$ is an isomorphism for any dual normed $\IR G$-module $V^\#$;
		\item The relative Johnson class $[J_\F]\in H^1_\Fb(G;(\linf(G/\F)/\IR)^\#)$ lies in the image of the canonical map $\can^1_{\F\subset\G,b}$.
	\end{enumerate}
	\begin{proof}
		Suppose that any subgroup $H\in \G$ is amenable relative to $\F|_H$. Then the resolution of normed $\IR G$-modules $V^\#\to C^*_\G(G;V^\#)$ is $\F$-strong and relatively $\F$-injective by Lemmas~\ref{lem:ex rel inj}~\ref{item:F-strong} and~\ref{lem:amen rel inj} applied to the $G$-set $G/\G$. Hence the canonical map $\can^n_{\F\subset\G,b}$ is an isomorphism for all $n\ge 0$ by Proposition~\ref{prop:fund lemma II}.
		
		The implications (ii) $\Rightarrow$ (iii) $\Rightarrow$ (iv) are obvious.
		Suppose that the relative Johnson class $[J_\F]$ lies in the image of the canonical map $\can^1_{\F\subset\G,b}$ and denote $V\coloneqq \linf(G/\F)/\IR$. We claim that for any subgroup $H\in \G$, the image of $[J_\F]$ under the restriction map 
		\[
			\res^1_{H\subset G,b}\colon H^1_\Fb(G;V^\#)\to H^1_{\F|_H,b}(H;V^\#)
		\]
		is trivial. Indeed, there is a commutative diagram
		\[\begin{tikzcd}
			H^1_{G,b}(\EGG;V^\#)\ar{rr}{\can^1_{\F\subset\G,b}}\ar{d} && H^1_{G,b}(\EFG;V^\#)\ar{d}\ar{dr}{\res^1_{H\subset G,b}} \\
			H^1_{H,b}(\EGG;V^\#)\ar{rr} && H^1_{H,b}(\EFG;V^\#)\ar{r}{\cong} & H^1_{H,b}(E_{\F|_H}H;V^\#) \,,
		\end{tikzcd}\]
		where the vertical maps are induced by viewing a $G$-space as an $H$-space. Observe that the lower left corner $H^1_{H,b}(\EGG;V^\#)$ is trivial, since when viewed as an $H$-space $\EGG$ is a model for $E_{\ALL|_H}H$ and hence $H$-equivariantly contractible. This proves the claim.
		
		Now, fix a subgroup $H\in\G$ and denote $W\coloneqq \linf(H/\F|_H)/\IR$. Consider the commutative diagram of normed $\IR H$-modules
		\[\begin{tikzcd}
			0\ar{r} & V^\#\ar{r}\ar{d} & \linf(G/\F)^\#\ar{r}\ar{d} & \IR\ar{r}\ar{d}{=} & 0 \\
			0\ar{r} & W^\#\ar{r} & \linf(H/\F|_H)^\#\ar{r} & \IR\ar{r} & 0 \,,
		\end{tikzcd}\]
		where the rows are exact, and remain exact when restricted to $L$-fixed-points for any $L\in\F|_H$. By Lemma~\ref{lem:basic properties} there are associated long exact sequences on bounded cohomology
		\[\begin{tikzcd}
			0\ar{r} & (V^\#)^H\ar{r}\ar{d} & (\linf(G/\F)^\#)^H\ar{r}\ar{d} & \IR\ar{r}{\partial^0_{V^\#}}\ar{d}{=} & H^1_{\F|_H,b}(H;V^\#)\ar{r}\ar{d} & \cdots \\
			0\ar{r} & (W^\#)^H\ar{r} & (\linf(H/\F|_H)^\#)^H\ar{r} & \IR\ar{r}{\partial^0_{W^\#}} & H^1_{\F|_H,b}(H;W^\#)\ar{r} & \cdots \,.
		\end{tikzcd}\]
		Observe that the image of $\partial^0_{V^\#}$ is precisely $\IR\cdot \res^1_{H\subset G,b}[J_\F]$ and hence trivial by the claim above. This implies that the map $\partial^0_{W^\#}$ is trivial and hence there exists a non-trivial $H$-invariant element in $\linf(H/\F|_H)^\#$. Thus $H$ is amenable relative to $\F|_H$ by Lemma~\ref{lem:rel amen equiv}. This finishes the proof.
	\end{proof}
\end{thm}

As special cases of Theorem~\ref{thm:unification} we obtain Theorem~\ref{thm:thm A} by taking $\G=\ALL$ and Theorem~\ref{thm:thm C} by taking $\F=\TR$. The case when $\F=\TR$ and $\G=\ALL$ recovers Theorem~\ref{thm:Johnson}.

\begin{cor}\label{cor:Loeh-Sauer}
	Let $X$ be a CW-complex with fundamental group $G$ and $\F$ be a family consisting of amenable subgroups of $G$. Suppose that there exists a model for $\EFG$ whose orbit space $G\backslash\EFG$ is homotopy equivalent to a $k$-dimensional CW-complex. Then the comparison map $c^n_X\colon H^n_b(X;\IR)\to H^n(X;\IR)$ vanishes for all $n>k$.
	\begin{proof}
		By Gromov's Mapping Theorem (see e.g.~\cite[Theorem 5.9]{Frigerio17}), the comparison map $c^n_X$ vanishes if the comparison map $c^n_{EG}\colon H^n_\Gb(EG;\IR)\to H^n_G(EG;\IR)$ vanishes. The $G$-map $EG\to \EFG$ induces a commutative square
		\[\begin{tikzcd}
			H^n_\Gb(EG;\IR)\ar{r}{c^n_{EG}} & H^n_G(EG;\IR) \\
	H^n_\Gb(\EFG;\IR)\ar{u}{\can^n_{\TR\subset\F,b}}[swap]{\cong}\ar{r}{c^n_{\EFG}} & H^n_G(\EFG;\IR)\ar{u}[swap]{\can^n_{\TR\subset\F}} \,,
		\end{tikzcd}\]
		where the canonical map $\can^n_{\TR\subset\F,b}$ is an isomorphism by Theorem~\ref{thm:thm C}. Since we are considering trivial coefficients, the lower right corner can be identified with the (non-equivariant) cohomology of the orbit space
		\[
			H^n_G(\EFG;\IR)\cong H^n(G\backslash\EFG;\IR) 
		\]
		(see e.g.~\cite[Theorem 4.2]{Fluch10}).
	\end{proof}
\end{cor}

As an application of Corollary~\ref{cor:Loeh-Sauer} we obtain the following well-known examples.

\begin{ex}\label{ex:factorization}
	The comparison map vanishes in all positive degrees for CW-complexes with the following fundamental groups:
	\begin{enumerate}[label=\enum]
		\item Graph products of amenable groups (e.g.\ right-angled Artin groups);
		\item Fundamental groups of graphs of groups with amenable vertex groups.
	\end{enumerate} 		
		Indeed, if $G_\Gamma$ is a graph product of amenable groups, we consider the family $\F$ generated by the vertex groups and direct products of vertex groups whenever the corresponding vertices form a clique in the underlying graph $\Gamma$. Then there exists a model for $E_\F (G_\Gamma)$ with contractible orbit space, which can be constructed by induction on the number of vertices of $\Gamma$. 
		
		If $G$ is the fundamental group of a graph of groups with amenable vertex groups, we consider the family $\F$ generated by the vertex groups. Then the Bass--Serre tree is a 1-dimensional model for $\EFG$. Recall that the comparison map always vanishes in degree 1, since $H^1_b(G;\IR)$ is trivial for any group $G$.
\end{ex}

We also obtain a characterization of relative amenability via relatively $\F$-injective modules, analogous to~\cite[Proposition 4.18]{Frigerio17}.

\begin{prop}\label{prop:rel amen rel inj}
	Let $G$ be a group and $\F$ be a family of subgroups. The following are equivalent:
	\begin{enumerate}[label=\enum]
		\item $G$ is amenable relative to $\F$;
		\item Any dual normed $\IR G$-module $V^\#$ is relatively $\F$-injective;
		\item The trivial normed $\IR G$-module $\IR$ is relatively $\F$-injective.
	\end{enumerate}
	\begin{proof}
		Suppose that $G$ is amenable relative to $\F$ and let $m_\F$ be a $G$-invariant mean on $G/\F$. The inclusion $V^\#\to C^0_\Fb(G;V^\#)$ of normed $G$-modules admits a right inverse $r$ given by
		\[
			r(f)(v) = m_\F\big(gH\mapsto f(gH)(v)\big)
		\]
		for $f\in C^0_\Fb(G;V^\#)$ and $v\in V$. Then the relative $\F$-injectivity of $V^\#$ follows from the relative $\F$-injectivity of $C^0_\Fb(G;V^\#)$.
		
		Clearly, condition (ii) implies (iii). Suppose that $\IR$ is relatively $\F$-injective. Consider the strongly $\F$-injective map $i\colon \IR\to \linf(G/\F)$ of normed $G$-modules that has an $H$-section $\tau_H$ given by $\tau_H(f)=f(eH)$ for each $H\in\F$. Then the identity $\id_\IR$ admits an extension along $i$ which yields a non-trivial $G$-invariant element in $\linf(G/\F)^\#$. By Lemma~\ref{lem:rel amen equiv} this finishes the proof.
	\end{proof}
\end{prop}

\paragraph{Characterization of relative finiteness}
Analogous to Theorem~\ref{thm:unification}, when instead considering all (not necessarily dual) normed $\IR G$-modules, one proves the theorem below. Let $G$ be a group and $\F$ be a family of subgroups.

Let $\ell^1(G/\F)$ denote the normed $\IR G$-module of summable functions $f\colon G/\F\to \IR$ with norm $\|f\|_1=\sum_{gH\in G/\F}f(gH)$ and let $\ell^1_0(G/\F)$ be the kernel of $\|\cdot\|_1$. We define the class $[K_\F]\in H^1_\Fb(G;\ell^1_0(G/\F))$ as the cohomology class of the 1-cocycle $K_\F\in C^1_\Fb(G;\ell^1_0(G/\F))$ given by 
\[
	K_\F(g_0H_0,g_1H_1)=\chi_{g_1H_1}-\chi_{g_0H_0} \,,
\]
where $\chi_{g_iH_i}$ is the characteristic function supported at $g_iH_i$ for $i=0,1$.

We say that $G$ is \emph{finite relative} to $\F$, if $\F$ contains a finite index subgroup of $G$.

\begin{thm}\label{thm:unification finite}
	Let $G$ be a group and $\F\subset\G$ be two families of subgroups. The following are equivalent:
	\begin{enumerate}[label=\enum]
		\item Any subgroup $H\in\G$ is finite relative to $\F|_H$;
		\item The canonical map $H^n_{\G,b}(G;V)\to H^n_\Fb(G;V)$ is an isomorphism for any normed $\IR G$-module $V$ and all $n\ge 0$;
		\item The canonical map $H^1_{\G,b}(G;V)\to H^1_\Fb(G;V)$ is an isomorphism for any normed $\IR G$-module $V$;
		\item The class $[K_\F]\in H^1_\Fb(G;\ell^1_0(G/\F))$ lies in the image of the canonical map $\can^1_{\F\subset\G,b}$.
	\end{enumerate}
\end{thm}
Theorem~\ref{thm:unification finite} has the following interesting special cases.
If $\F$ is arbitrary and $\G=\ALL$, we characterize that $\F$ contains a finite index subgroups of $G$. If $\F=\TR$ and $\G$ is arbitrary, we characterize that all subgroups in $\G$ are finite, generalizing~\cite[Theorem B]{Moraschini-Raptis21}. We recover the characterization of finite groups (\cite[Theorem 3.12]{Frigerio17}) for $\F=\TR$ and $\G=\ALL$.

\section{Characterization of relative hyperbolicity}
In this section we prove a characterization of relatively hyperbolic groups in terms of bounded Bredon cohomology analogous to Theorem~\ref{thm:Mineyev}.

Let $G$ be a finitely generated group and $\calH$ be a finite set of subgroups. Recall that $G$ is \emph{hyperbolic relative to $\calH$} if the coned--off Cayley graph is hyperbolic and fine. For example, hyperbolic groups are hyperbolic relative to the trivial subgroup, free products $G_1\ast G_2$ are hyperbolic relative to $\{G_1,G_2\}$, and fundamental groups of finite volume hyperbolic manifolds are hyperbolic relative to the cusp subgroups. If $G$ is hyperbolic relative to $\calH$, then $\calH$ is almost malnormal and hence malnormal if $G$ is torsionfree.

From now on, let the ring $R$ be either $\IQ$ or $\IR$. A map $f\colon C\to B$ of normed $RG$-modules is called \emph{undistorted} if there exists a constant $K\ge 0$ such that for all $b\in \im(f)$ there exists $c\in C$ with $f(c)=b$ such that $\|c\|_C\le K\cdot \|b\|_B$.
A normed $RG$-module $P$ is called \emph{boundedly projective} if for every undistorted epimorphism $f\colon C\to B$ and every map $\phi\colon P\to B$ of normed $RG$-modules, there exists a map $\Phi\colon P\to C$ of normed $RG$-modules such that $f\circ \Phi=\phi$.

The following lemma~\cite[Lemma 52]{Mineyev-Yaman07} is useful to construct $G$-equivariant maps.

\begin{lem}[Mineyev--Yaman]\label{lem:b-projective}
	Let $G$ be a group and $S$ be a $G$-set with finite stabilizers. Then $\IQ[S]$ is projective as a $\IQ G$-module and boundedly projective as a normed $\IQ G$-module when equipped with the $\ell^1$-norm. 
\end{lem}

Let $X$ be a $G$-CW-complex with cocompact $(n+1)$-skeleton and consider for $k\ge 0$ the cellular chains $C^\cell_k(X;R)$ as a normed $RG$-module equipped with the $\ell^1$-norm. We say that $X$ satisfies a \emph{linear homological isoperimetric inequality over $R$ in degree $n$} if the boundary map
\[
	\partial_{n+1}\colon C^\cell_{n+1}(X;R)\to C^\cell_n(X;R)
\]
is undistorted.
Equivalently, there exists a constant $K\ge 0$ such that for every cellular $n$-boundary $b\in B^\cell_n(X;R)$ we have $\|b\|_\partial\le K\cdot \|b\|_1$, where
\[
	\|b\|_\partial\coloneqq \inf\{\|c\|_1\ |\ c\in C^\cell_{n+1}(X;R),\, \partial_{n+1}(c)=b\}
\]
(which is sometimes called the \emph{filling norm}).

If $G$ is hyperbolic relative to $\calH$, Mineyev--Yaman~\cite[Theorem 41]{Mineyev-Yaman07} have constructed the so called ``ideal complex" $X$. It is in particular a cocompact $G$-CW-complex with precisely one equivariant 0-cell $G/H$ for each each $H\in\calH$ and finite edge-stabilizers. Moreover, $X$ is (non-equivariantly) contractible and hence a model for $E_\FcalH G$ provided that $G$ is torsionfree. 
We summarize some of its properties~\cite[Theorem 47 and 51]{Mineyev-Yaman07}.

\begin{thm}[Mineyev--Yaman]\label{thm:ideal complex}
	Let $G$ be a finitely generated torsionfree group and $\calH$ be a finite set of subgroups. If $G$ is hyperbolic relative to $\calH$, then there exists a cocompact model $X$ for $E_\FcalH G$ satisfying the following:
	\begin{enumerate}[label=\enum]
		\item\label{item:ideal linear} $X$ satisfies linear homological isoperimetric inequalities over $\IQ$ in degree $n$ for all $n\ge 1$;
		\item\label{item:ideal bicombing} There exists a map $q\colon X^{(0)}\times X^{(0)}\to C_1^\cell(X;\IQ)$ with $\partial_1(q(a,b))=b-a$, called a homological $\IQ$-bicombing, that is $G$-equivariant and satisfies
		\[
			\|q(a,b)+q(b,c)-q(a,c)\|_1\le K
		\]
		for all $a,b,c\in X^{(0)}$ and a uniform constant $K\ge 0$.
	\end{enumerate}
\end{thm}

The following criterion for relative hyperbolicity is a combination of~\cite[Proposition 8.3 and Theorem 8.5]{Franceschini18} (see also~\cite[Theorems 1.6 and 1.10]{MartinezPedroza16}).

\begin{thm}[Franceschini, Mart\'inez-Pedroza]\label{thm:Franceschini}
	Let $G$ be a group and $\calH$ be a finite set of subgroups. Then $G$ is hyperbolic relative to $\calH$ if there exists a $G$-CW-complex $Z$ satisfying the following:
	\begin{enumerate}[label=\enum]
		\item $Z$ is simply-connected;
		\item The 2-skeleton $Z^{(2)}$ is cocompact;
		\item $\calH$ is a set of representatives of distinct conjugacy classes of vertex-stabilizers such that each infinite stabilizer is represented;
		\item The edge-stabilizers of $Z$ are finite;
		\item $Z$ satisfies a linear homological isoperimetric inequality over $\IR$ in degree 1.
	\end{enumerate}
\end{thm}

We prove the following characterization of relative hyperbolicity closely following Mineyev's original proof of Theorem~\ref{thm:Mineyev} (\cite[Theorem 11]{Mineyev01} and~\cite[Theorem 9]{Mineyev02}).

\begin{thm}\label{thm:rel hyp}
	Let $G$ be a finitely generated torsionfree group and $\calH$ be a finite malnormal collection of subgroups. Let $\F$ be the family $\FcalH$ and suppose that $G$ is of type $F_{2,\F}$. Then the following are equivalent:
	\begin{enumerate}[label=\enum]
		\item\label{item:rel hyp i} $G$ is hyperbolic relative to $\calH$;
		\item\label{item:rel hyp ii} The comparison map $H^n_\Fb(G;V)\to H^n_\F(G;V)$ is surjective for any normed $\IQ G$-module $V$ and all $n\ge 2$;
		\item\label{item:rel hyp iii} The comparison map $H^2_\Fb(G;V)\to H^2_\F(G;V)$ is surjective for any normed $\IR G$-module $V$.
	\end{enumerate}
	\begin{proof}
		Suppose that $G$ is hyperbolic relative to $\calH$.
		Let $X$ be the model for $\EFG$ that is given by Mineyev--Yaman's ideal complex (Theorem~\ref{thm:ideal complex}) and $Y$ be the simplicial model for $\EFG$ with (non-equivariant) $n$-cells corresponding to $(G/\F)^{n+1}$ for all $n\ge 0$.
		We claim that there is a $G$-chain map
		\[
			\varphi_*\colon C^\cell_*(Y;\IQ)\to C^\cell_*(X;\IQ)
		\]
		with $\varphi_n$ bounded for all $n\ge 2$, admitting a $G$-homotopy left inverse. We construct $\varphi_*$ inductively as follows. In degree 0, we define
		\[
			\varphi_0\colon C^\cell_0(Y;\IQ)=\IQ[G/\F]\to C^\cell_0(X;\IQ)
		\]
		to map a generator of the form $eH$ to the vertex of $X$ with stabilizer containing $H$. Then extend $G$-equivariantly and $\IQ$-linearly to all of $\IQ[G/\F]$.
		In degree 1, we define $\varphi_1\colon C^\cell_1(Y;\IQ)\to C^\cell_1(X;\IQ)$ on generators by
		\[
			\varphi_1(g_0H_0,g_1H_1)=q(\varphi_0(g_0H_0),\varphi_0(g_1H_1)) \,,
		\]
		where $q$ is the homological $\IQ$-bicombing on $X$ from Theorem~\ref{thm:ideal complex}~\ref{item:ideal bicombing}. Since both $q$ and $\varphi_0$ are $G$-equivariant, so is $\varphi_1$. In degree 2, we consider the maps
		\begin{equation}\label{eqn:q2}\begin{tikzcd}
			C^\cell_2(Y;\IQ)\ar{r}{\partial_2^Y} & C^\cell_1(Y;\IQ)\ar{d}{\varphi_1} \\
			C^\cell_2(X;\IQ)\ar{r}{\partial_2^X} & C^\cell_1(X;\IQ)
		\end{tikzcd}\end{equation}
		and observe that the composition $\varphi_1\circ \partial_2^Y$ is bounded by properties of $q$ and that $\partial_2^X$ is undistorted by Theorem~\ref{thm:ideal complex}~\ref{item:ideal linear}. There is a $G$-invariant decomposition
		\[
			C^\cell_2(Y;\IQ) \cong \IQ[S_1]\oplus \IQ[S_2] \,,
		\]
		where $S_1$ and $S_2$ denote the sets of 2-cells of $Y$ with trivial resp.\ non-trivial stabilizers.
		We obtain a bounded $G$-map $\varphi_2\colon C^\cell_2(Y;\IQ)\to C^\cell_2(X;\IQ)$ by using the bounded projectivity of $\IQ[S_1]$ (Lemma~\ref{lem:b-projective}) and by setting $\varphi_2$ to be zero on $\IQ[S_2]$.
		This renders the square~\eqref{eqn:q2} commutative because the edge-stabilizers of $X$ are trivial.
		
		Assuming that $\varphi_n$ has been constructed, one analogously defines a bounded $G$-map $\varphi_{n+1}$ using that $\partial_{n+1}^X$ is undistorted by Theorem~\ref{thm:ideal complex}~\ref{item:ideal linear}. Thus one obtains a $G$-chain map $\varphi_*$ with $\varphi_n$ bounded for $n\ge 2$. To conclude the claim, we note that $C^\cell_*(Y;\IQ)$ is a relatively $\F$-projective $\F$-strong resolutions of $\IQ$ by Lemma~\ref{lem:rel proj}. Hence by Proposition~\ref{prop:fund lemma I} any $G$-chain map $\psi_*\colon C^\cell_*(X;\IQ)\to C^\cell_*(Y;\IQ)$ extending $\id_\IQ$ is a $G$-homotopy left inverse of $\varphi_*$.
		
		Now, let $V$ be a normed $\IQ G$-module. Applying $\Hom_{\IQ G}(-,V)$ yields a cochain map
		\[
			\varphi^*\colon C^*_\cell(X;V)^G\to C^*_\cell(Y;V)^G 
		\]
		with homotopy right inverse $\psi^*$. In particular, the composition $\varphi^*\circ \psi^*$ induces the identity on $H^*(C^*_\cell(Y;V)^G)\cong H^*_\F(G;V)$. Finally, for $n\ge 2$ let $c\in C^n_\cell(Y;V)^G$ be a cocycle. Then $\varphi^n(\psi^n(c))$ and $c$ represent the same cohomology class in $H^n_\F(G;V)$. We have
		\[
			\|\varphi^n(\psi^n(c))\|_\infty = \|\psi^n(c)\circ \varphi_n\|_\infty \le \|\psi^n(c)\|_\infty \cdot \|\varphi_n\|_\infty \,,
		\]
		where $\varphi_n$ is bounded by construction and so is $\psi^n(c)\in C^n_\cell(X;V)^G$ because $X$ has only finitely many orbits of $n$-cells. Thus we have shown that for $n\ge 2$ every cohomology class in $H^n_\F(G;V)$ admits a bounded representative.
	\\
	
	Obviously condition (ii) implies (iii).	Suppose that the comparison map is surjective in degree 2 for coefficients in any normed $\IR G$-module. Let $Z$ be a model for $\EFG$ with cocompact 2-skeleton. 
	Since $\calH$ is malnormal, by collapsing fixed-point sets of $Z$ we may assume that for every non-trivial subgroup $H\in\F$ the fixed-point set $Z^H$ consists of precisely one point.\footnote{A detailed proof of this fact was communicated to us by Sam Hughes. It will appear in a forthcoming paper of his, joint with Mart\'inez-Pedroza and S\'{a}nchez Salda\~{n}a.}
	In order to apply Theorem~\ref{thm:Franceschini} and conclude that $G$ is hyperbolic relative to $\calH$, it remains to verify that $Z$ satisfies a linear homological isoperimetric inequality over $\IR$ in degree 1.
	
	We take as coefficients the cellular 1-boundaries $V\coloneqq B^\cell_1(Z;\IR)$ equipped with the norm $\|\cdot\|_\partial$. Let $Y$ be the simplicial model for $\EFG$. Then there is a $G$-chain homotopy equivalence
	\[
		\psi_*\colon C^\cell_*(Z;\IR)\to C^\cell_*(Y;\IR)
	\]
	with $G$-homotopy inverse $\varphi_*$. Applying $\Hom_{\IR G}(-,V)$ yields a cochain homotopy equivalence
	\[
		\psi^*\colon C^*_\cell(Y;V)^G\to C^*_\cell(Z;V)^G
	\]
	with homotopy inverse $\varphi^*$. In particular, the composition $\psi^*\circ \varphi^*$ induces the identity on $H^*(C^*_\cell(Z;V)^G)\cong H^*_\F(G;V)$.
	Consider the 2-cocycle $u\in C^2_\cell(Z;V)^G$ given by the boundary map
	\[
		u=\partial_2\colon C_2^\cell(Z;\IR)\to B_1^\cell(Z;\IR)=V \,.
	\]
	Then we can write
	\begin{equation}\label{eqn:u}
		u=(\psi^2\circ\varphi^2)(u)+\delta^1_Z(v)
	\end{equation}
	for some $v\in C^1_\cell(Z;V)^G$. Since the comparison map $H^2_\Fb(G;V)\to H^2_\F(G;V)$ is surjective by hypothesis, we can write
	\begin{equation}\label{eqn:u'}
		\varphi^2(u)=u'+\delta^1_Y(v')
	\end{equation}
	for a \emph{bounded} 2-cocycle $u'\in C^2_\cell(Y;V)^G$ and some $v'\in C^1_\cell(Y;V)^G$.
	For a fixed vertex $y\in Y^{(0)}=G/\F$, let $\Cone_{y}\colon C_1^\cell(Y;\IR)\to C_2^\cell(Y;\IR)$ be defined on generators by
	\[
		\Cone_{y}\left((g_0H_0,g_1H_1)\right)=(y,g_0H_0,g_1,H_1) \,.
	\]
	Obviously $\Cone_{y}$ preserves the $\ell^1$-norms.
	For a $G$-CW-complex $W$, we denote the evaluation pairing by 
	\[
		\spann{\cdot,\cdot}_W\colon C^*_\cell(W;V)^G\times C_*^\cell(W;\IR)\to V \,.
	\]
	Now, for $b\in C_1^\cell(Z;\IR)$ and $c\in C_2^\cell(Z;\IR)$ with $\partial_2(c)=b$, we find by~\eqref{eqn:u} that
	\begin{align*}
		b &= \spann{u,c}_Z = \spann{(\psi^2\circ\varphi^2)(u)+\delta^1_Z(v),c}_Z \\
		&= \spann{(\psi^2\circ\varphi^2)(u),c}_Z + \spann{v,\partial_2^Y(c)}_Z \\
		&= \spann{\varphi^2(u),\psi_2(c)}_Y + \spann{v,b}_Z \,.
	\end{align*}	
	Since $\varphi^2(u)$ is a cocycle and $\psi_2(c)-\Cone_y(\partial_2^Y(\psi_2(c)))$ is a cycle and hence a boundary, we have
	\begin{align*}
		\spann{\varphi^2(u),\psi_2(c)}_Y
		&= \spann{\varphi^2(u),\Cone_y(\partial_2^Y(\psi_2(c)))}_Y \\
		&= \spann{\varphi^2(u),\Cone_y(\psi_1(b))}_Y \\
		&= \spann{u'+\delta^1_Y(v'),\Cone_y(\psi_1(b))}_Y \\
		&= \spann{u',\Cone_{y}(\psi_1(b))}_Y + \spann{v',\partial_2^Y(\Cone_{y}(\psi_1(b)))}_Y\\
		&= \spann{u',\Cone_{y}(\psi_1(b))}_Y + \spann{v',\psi_1(b)}_Y \\
		&= \spann{u',\Cone_{y}(\psi_1(b))}_Y + \spann{\psi^1(v'),b}_Z \,,
	\end{align*}
	where we used~\eqref{eqn:u'}.
	Together, we have
	\[
		b=\spann{u',\Cone_y(\psi_1(b))}_Y + \spann{\psi^1(v')+v,b}_Z
	\]
	and it follows that
	\begin{align*}
		\|b\|_\partial &\le \|\spann{u',\Cone_y(\psi_1(b))}_Y\|_\partial + \|\spann{\psi^1(v')+v,b}_Z\|_\partial \\
		&\le \|u'\|_\infty\cdot \|\Cone_y(\psi_1(b))\|_1 + \|\psi^1(v')+v\|_\infty \cdot \|b\|_1 \\
		&= \|u'\|_\infty\cdot \|\psi_1(b)\|_1 + \|\psi^1(v')+v\|_\infty \cdot \|b\|_1 \\
		&\le (\|u'\|_\infty\cdot \|\psi_1\|_\infty+\|\psi^1(v')+v\|_\infty)\cdot \|b\|_1 \,.
	\end{align*}	
	Finally, $u'$ is bounded by construction and so are $\psi_1$ and $\psi^1(v')+v$ because they are $G$-maps with domain $C_1^\cell(Z;\IR)$ and $Z$ has only finitely many orbits of 1-cells. Thus we have shown that $Z$ satisfies a linear homological isoperimetric inequality over $\IR$ in degree 1. This finishes the proof.
	\end{proof}
\end{thm}

\begin{rem}[Groups with torsion]\label{rem:torsion}
	In Theorem~\ref{thm:rel hyp}, if the group $G$ is not assumed to be torsionfree and $\calH$ is instead assumed to be almost malnormal, one can still prove the equivalence of~\ref{item:rel hyp i} and~\ref{item:rel hyp iii}. However, a few modifications are necessary which we shall only outline. 
	
	Assuming that $G$ is hyperbolic relative to $\calH$, Mineyev--Yaman's ideal complex has to be replaced by a Rips type construction $X$ due to Mart\'inez-Pedroza--Przytycki that is a model for $E_{\F\cup\FIN}G$. This complex $X$ satisfies a linear homological isoperimetric inequality over $\IZ$ in degree 1 (\cite[Corollary 1.5]{MartinezPedroza19}). It is part of a hyperbolic tuple in the sense of~\cite[Definition 38]{Mineyev-Yaman07} and hence admits a homological $\IQ$-bicombing by~\cite[Theorem 47]{Mineyev-Yaman07}. Then one can construct a $G$-chain map $\varphi_*$ with $\varphi_2$ bounded similarly as before and conclude surjectivity of the comparison map in degree 2 for the family $\F\cup\FIN$. This implies the same for the family $\F$ over the ring $\IR$.
	
	For the converse implication, since $\calH$ is almost malnormal, there exists a model $Z$ for $\EFG$ such that for every infinite subgroup $H\in\F$ the fixed-point set $Z^H$ consists of precisely one point. Then one shows as before that $Z$ satisfies a linear homological isoperimetric inequality over $\IR$ in degree 1 and concludes by Theorem~\ref{thm:Franceschini}.
		
	We do not know whether condition~\ref{item:rel hyp ii} is equivalent to~\ref{item:rel hyp i} and~\ref{item:rel hyp iii} in this case.
\end{rem}

\bibliographystyle{alpha}
\bibliography{bib}

\end{document}